\documentclass[12pt,a4paper]{article}

\usepackage[latin2]{inputenc}
\usepackage{epsf}
\usepackage{amsmath}
\usepackage{amsthm}
\usepackage{amsfonts}
\usepackage{amssymb}
\usepackage{setspace}

\def\p{{\bf P}}

\def\eq#1{\smash{\mathop{=}\limits^{#1}}}

\begin{document}

\begin{center}
{\large\bf A Supplement to the Paper Poisson Approximation in a Poisson Limit Theorem Inspired by Coupon Collecting}
\end{center}
\smallskip
\begin{center}
Anna P\'{o}sfai\\
\smallskip
{\small \noindent Analysis and Stochastics Research Group of the Hungarian Academy of Sciences, Bolyai Institute, University of Szeged, Aradi v\'ertan\'uk tere 1, Szeged 6720, Hungary; e-mail: posfai@math.u-szeged.hu}
\end{center}

\bigskip\smallskip

In this note we give a proof for the result stated as Theorem 4 in \cite{P}.

A collector samples with replacement a set of $n\in N:=\{1,2,\ldots\}$ distinct coupons so that the draws are independent and at each time any one of the $n$ coupons is drawn with the same probability $1/n$. For an integer $m_n\in\{0,1,\ldots,n-1\}$ that depends on~$n$, sampling is repeated until the first
time $W_{n,m_n}$ that the collector has collected $n-m_n$ distinct coupons.
Baum and Billingsley proved in \cite{BB} (using the method of characteristic functions) that if
\begin{equation}\label{m_n}
m_n\to\infty \quad\textrm{ and }\quad\frac{n-m_n}{\sqrt{n}}\to\sqrt{2\lambda}\quad\textrm{ for some }\lambda>0\textrm{ constant, as } n\to\infty,
\end{equation}
then $W_{n,m_n}-(n-m_n)$ converges in distribution to the Poisson law with mean $\lambda$.

Throughout all asymptotic relations are meant as $n\to\infty$.

It can be seen that the following equality in distribution holds for $\widetilde{W}_{n,m_n}:=W_{n,m_n}-(n-m_n)$:
$$\widetilde{W}_{n,m_n}\;\eq{\cal D}\;\sum_{i=m_n+1}^{n}\widetilde{X}_{n,i}$$
where the $\widetilde{X}_{ni}$ random variables are independent, and $\widetilde{X}_{n,i}+1$ has geometric distributions with success probability $i/n$, $i\in\{m_n+1,\ldots, n\}$, $n\in N$, that is $\p\{\widetilde{X}_{n,i}+1 = j\} = \left(1 - \frac{i}{n}\right)^{j-1}\frac{i}{n}$, $j\in N$, $i\in\{m_n+1,\ldots, n\}$.

We approximate the waiting time $\widetilde{W}_{n,m_n}$ with a Poisson random variable that has mean $\lambda_n=\sum_{i=m_n+1}^{n}\left(1-\frac{i}{n}\right)$. Using the special combinatorial structure of the problem we derive the first asymptotic correction of the $\p(\widetilde{W}_{n,m_n}=k)$, $k=0,1,\ldots$, probabilities to the corresponding Poisson point probabilities. We note that in principal the method presented in the proof can be extended to determine higher order terms in the asymptotic expansion.

\bigskip

\noindent\textbf{Theorem.} If $\{m_n\}_{n\in N}$ is a sequence of nonnegative integers that satisfies (\ref{m_n}),
\begin{equation*}
\lambda_n=\sum_{i=m_n+1}^{n}\left(1-\frac{i}{n}\right) \quad\textrm{and}\quad \lambda_{n,2}=\sum_{i=m_n+1}^{n}\left(1-\frac{i}{n}\right)^2,
\end{equation*}
then
\begin{align*}
&\p(\widetilde{W}_{n,m_n}=0)=e^{-\lambda_n}
-e^{-\lambda_n}\frac{\lambda_{n,2}}{2}
+O\left(\frac{1}{n}\right),\\
&\p(\widetilde{W}_{n,m_n}=1)=e^{-\lambda_n}\lambda_n
-e^{-\lambda_n}\lambda_n\frac{\lambda_{n,2}}{2}
+O\left(\frac{1}{n}\right),\\
&\p(\widetilde{W}_{n,m_n}=k)
=e^{-\lambda_n}\frac{\lambda_n^k}{k!}
+e^{-\lambda_n}\left(\frac{\lambda_n^{k-2}}{(k-2)!}-\frac{\lambda_n^{k}}{k!}\right)\frac{\lambda_{n,2}}{2}
+O\left(\frac{1}{n}\right),\quad k\geq2.
\end{align*}

\bigskip

We note that $\lambda_{n,2}=\frac{(2\lambda_n)^{3/2}}{3\sqrt{n}} + O\left(\frac{1}{n}\right)$. Indeed,
\begin{align*}
\lambda_{n,2}
&=\sum_{i=m_n+1}^{n}\left(1-\frac{i}{n}\right)^2\\
&=n-m_n-\frac{2}{n}\left[\frac{n(n+1)}{2}-\frac{m_n(m_n+1)}{2}\right]+\frac{1}{n^2}\left[\frac{n(n+1)(2n+1)}{2}-\frac{m_n(m_n+1)(2m_n+1)}{2}\right]\\
&=\frac{(n-m_n)(n-m_n-1)(n-m_n-\frac{1}{2})}{3n^2}\\
&=\frac{(2\lambda_n)^{3/2}}{3\sqrt{n}}+\left(\frac{(n-m_n)(n-m_n-1)}{3n^2}\left[n-m_n-\frac{1}{2}-\sqrt{(n-m_n)(n-m_n-1)}\right]\right),
\end{align*}
where we used the fact that $\lambda_n=\frac{(n-m_n)(n-m_n-1)}{2n}$, and the second term in the formula above is $O\left(\frac{1}{n}\right)$ by (\ref{m_n}).

We shall need the following simple result for the proof of the theorem.

\bigskip

\noindent\textbf{Proposition.} If $\{m_n\}_{n\in N}$ is a sequence of integers that satisfies (\ref{m_n}), then
\begin{align}
&\lambda_n=\lambda_{n,1}:=\sum_{i=m_n+1}^{n}\left(1-\frac{i}{n}\right)\to\lambda,\quad\textrm{ and}\label{lambda_n}\\
&\lambda_{n,j}:=\sum_{i=m_n+1}^{n}\left(1-\frac{i}{n}\right)^j\leq\lambda_n\left(\frac{2\lambda_n}{n}\right)^{\frac{j-1}{2}},
\quad\textrm{and}\quad\lambda_{n,j}\to0,\quad j=2,3,\ldots\label{lambda_n,j}
\end{align}

\bigskip

\noindent\textbf{Proof.} (\ref{lambda_n}) is true, because
\begin{equation*}
\lambda_n=\sum_{i=m_n+1}^{n}\left(1-\frac{i}{n}\right)
=n-m_n-\frac{1}{n}\left[\frac{n(n+1)}{2}-\frac{m_n(m_n+1)}{2}\right]
=\frac{(n-m_n)(n-m_n-1)}{2n}\to\lambda
\end{equation*}
by (\ref{m_n}). By taking the square root of both sides of the equality above it can be deduced that
\begin{equation}\label{csillag}
\frac{n-m_n-1}{\sqrt{n}}\leq\sqrt{2\lambda_n}.
\end{equation}
Now we prove the first assertion of (\ref{lambda_n,j}) by induction. For an arbitrary $j=2,3,\ldots$ we bound $\lambda_{n,j}$ as follows:
$$\lambda_{n,j}=\sum_{i=m_n+1}^{n}\left(1-\frac{i}{n}\right)^j
\leq\frac{n-m_n-1}{n}\sum_{i=m_n+1}^{n}\left(1-\frac{i}{n}\right)^{j-1}
=\frac{n-m_n-1}{n}\lambda_{n,j-1}
$$
Since for $j=2$ this gives $\lambda_{n,2}\leq\lambda_n\sqrt{\frac{2\lambda_n}{n}}$ by (\ref{csillag}), we have the first part of (\ref{lambda_n,j}) in this case. If we have the same result for some $j>2$, then it holds true for $j+1$ as well by the argument above, (\ref{csillag}) and the inductional hypothesis. Since $\lambda_n\to\lambda$ by (\ref{lambda_n}), the second part of (\ref{lambda_n,j}) follows from the first. $\Box$

\bigskip

\noindent\textbf{Proof of the Theorem.} We are going to represent each possible outcome of the collector's sampling with a sequence of integers the following way: let us suppose that while sampling (with replacement), the collector labels the distinct coupons he draws form 1 to $n-m_n$ in the order he obtains them in the course of time, and after each draw he writes down the label of the coupon just drawn. So he begins the enumeration of labels with a 1 after the first draw, and each number that he writes to the end of his list after a draw is either the label already on the coupon he just got (if he had drawn the same one before), or it is the label he gives the coupon at that moment, which would be the smallest positive integer he has not yet used in the process of sampling and labeling. In the first case we call the new member of the sequence "superfluous", while in the second case we call it a "first appearance".

We fix an arbitrary $k\in N$, and we suppose that $n$ so big that $n-m_n>k$ holds. Now $\widetilde{W}_{n,m_n}=k$ means that the collector had $k$ "superfluous" draws, thus the corresponding representing sequence contains $n-m_n$ "first appearances" and $k$ "superfluous" members. We categorize all such outcomes according to how the $k$ "superfluous" draws are split into blocks by the $n-m_n$ "first appearances" in the representing sequences: to each vector $\underline{k}=(k_{m_n+1},k_{m_n+2},\ldots,k_{n-1})$, where $k_i\in Z_+$, $i=m_n+1,\ldots,n-1$, and $\sum_{i=m_n+1}^{n-1}k_i=k$, correspond the sequences where there are $k_{n-1}$ "superfluous" members between the 1st and 2nd "first appearances", $k_{n-2}$ "superfluous" members between the 2nd and 3rd "first appearances", and so on, $k_{m_n+1}$ "superfluous" members between the $(n-m_n-1)$th and $(n-m_n)$th "first appearances". (This is the same as saying that $\widetilde{X}_{ni}=k_i$, for all $i=m_n+1,\ldots,n$.) The probability of getting such a sequence is
\begin{multline*}
\frac{n}{n}\left(1-\frac{n-1}{n}\right)^{k_{n-1}}
\frac{n-1}{n}\left(1-\frac{n-2}{n}\right)^{k_{n-2}}\cdots
\left(1-\frac{m_n+1}{n}\right)^{k_{m_n+1}}\frac{m_n+1}{n}\\
=\left(\prod_{i=m_n+1}^{n}\frac{i}{n}\right)\prod_{i=m_n+1}^{n-1}\left(1-\frac{i}{n}\right)^{k_i}.
\end{multline*}
It follows that
\begin{equation}\label{p(w=k)}
\p(\widetilde{W}_{n,m_n}=k)=\left(\prod_{i=m_n+1}^{n}\frac{i}{n}\right)\sum_{\underline{k}\in I_k}\prod_{i=m_n+1}^{n-1}\left(1-\frac{i}{n}\right)^{k_i},
\end{equation}
where
$$I_k:=\left\{\underline{k}\in Z_+^{n-m_n-1} : \sum_{i=m_n+1}^{n-1}k_i=k\right\}.$$

\bigskip

Now we are going to examine the sum in (\ref{p(w=k)}) above, which we denote by $S_{n,m_n,k}=S_k$. For $k=0$ it is an empty sum, and thus it equals 1 by definition. Now let us suppose that $k>2$, we are going to return to the cases $k=0$ and 1 later on. For an arbitrary such $k$ we see that
$$I_k=\cup_{l=1}^{k} I_{k,l},\quad\textrm{ where }\quad
I_{k,l}=\{\underline{k}\in I_k : \underline{k} \textrm{ has exactly } l \textrm{ nonzero components}\}, l=1,\ldots,k,$$
and we correspondingly define $S_{k,l}$ to be the part of $S_{k}$ that contains the summands over $\underline{k}\in I_{k,l}$, thus we have
\begin{equation}\label{s_nm}
S_{k}=\sum_{\underline{k}\in I_k}\prod_{i=m_n+1}^{n-1}\left(1-\frac{i}{n}\right)^{k_i}
=\sum_{l=1}^{k}\sum_{\underline{k}\in I_{k,l}}\prod_{i=m_n+1}^{n-1}\left(1-\frac{i}{n}\right)^{k_i}
=\sum_{l=1}^{k}S_{k,l}.
\end{equation}

To determine the limit of $S_k$ we examine the asymptotic behavior of the $S_{k,l}$ expressions separately. We fix an arbitrary $l=1,\ldots,k$, and with $|A|$ denoting the cardinality of an arbitrary set $A$, we now calculate $|I_{k,l}|$. We can think of the vectors in $I_k$ as the results of distributing $k$ 1-s in $n-m_n-1$ spaces in all possible ways: to each of these distributions correspond a vector in $I_k$ whose $i$th component is the number of 1-s put in the $i$th space, $i=m_n+1,\ldots,n$. To produce a vector in $I_{k,l}$ we first choose $l$ different spaces, and we put a 1 in each of them, then we distribute the remaining $k-l$ 1-s in these previously chosen $l$ spaces that already have a 1, but this time any such space can be chosen more than once. This gives
$$|I_{k,l}|={n-m_n-1\choose l}{k-1\choose k-l},\quad l=1,\ldots,k.$$
We obviously bound $S_{k,l}$ from above if we replace each of the factors in its products by the largest one of them, namely by $1-\frac{m_n+1}{n}$. This together with the just calculated formula gives
\begin{equation*}
S_{k,l}
\leq{n-m_n-1\choose l}{k-1\choose k-l}\left(1-\frac{m_n+1}{n}\right)^k\\
\leq(k-1)!\left(\frac{n-m_n-1}{\sqrt{n}}\right)^{k+l}\left(\frac{1}{\sqrt{n}}\right)^{k-l}.
\end{equation*}
Hence by (\ref{csillag}) we have
\begin{equation}\label{s_maradek}
S_{k,l}\leq\frac{(k-1)!}{l!(l-1)!}\sqrt{2\lambda_n}^{k+l}\left(\frac{1}{\sqrt{n}}\right)^{k-l}\quad\textrm{and}\quad
\sum_{l=1}^{l'}S_{k,l}\leq k!\min\left\{1,(2\lambda_n)^{k}\right\}\left(\frac{1}{\sqrt{n}}\right)^{k-l'}
\end{equation}
for any $l'\in\{1,\ldots,k\}$. We see from the first inequality that $S_{k,l}$ goes to 0 for $l=1,\ldots,k-1$, but it gives a constant upper bound for $l=k$. We are going to examine the latter case more carefully. Notice that the components of a vector in $I_{k,k}$ are all 0-s and 1-s, thus for any $\underline{k}\in I_{k,k}$ $\frac{1}{k_{m_n+1}!k_{m_n+2}!\dots k_{n-1}!}=1$. Using this and the decomposition of the index set $I_k=\cup_{l=1}^{k}I_{k,l}$ we obtain
\begin{multline*}
S_{k,k}=\frac{1}{k!}\sum_{\underline{k}\in I_k}\frac{k!}{k_{m_n+1}!k_{m_n+2}!\dots k_{n-1}!}\prod_{i=m_n+1}^{n-1}\left(1-\frac{i}{n}\right)^{k_i}-\\
-\sum_{l=1}^{k-1}\sum_{\underline{k}\in I_{k,l}}\frac{1}{k_{m_n+1}!k_{m_n+2}!\dots k_{n-1}!}\prod_{i=m_n+1}^{n-1}\left(1-\frac{i}{n}\right)^{k_i}.
\end{multline*}
The first term of $S_{k,k}$ is equal to $\frac{1}{k!}\left[\sum_{i=m_n+1}^{n}\left(1-\frac{i}{n}\right)\right]^k$ by the polynomial theorem, thus we have
\begin{equation}\label{s_fotag}
S_{k,k}=\frac{\lambda_n^k}{k!}
-\sum_{l=1}^{k-1}\sum_{\underline{k}\in I_{k,l}}\frac{1}{k_{m_n+1}!k_{m_n+2}!\dots k_{n-1}!}\prod_{i=m_n+1}^{n-1}\left(1-\frac{i}{n}\right)^{k_i}.
\end{equation}
It follows that $\lim_{n\to\infty}S_{k,k}=\frac{\lambda^k}{k!}$, because we have (\ref{lambda_n}), and the sum above can be bounded by $\sum_{l=1}^{k-1}S_{k,l}$, which goes to 0 by (\ref{s_maradek}). Thus putting together our results for the expressions $S_{k,l}$ in (\ref{s_nm}), we conclude that the part of $S_k$ that counts -- in the sense that it asymptotically contributes a positive constant to $S_k$ --, is $S_{k,k}$, which is the part of the sum in the defining formula of $S_k$ that corresponds to the 0 - 1 vectors of the $I_k$ index set.

If we write (\ref{s_fotag}) into (\ref{s_nm}), we obtain the following formula for $S_k$:
\begin{equation}\label{s_k}
S_{k}=\frac{\lambda_n^k}{k!}+\sum_{l=1}^{k-1}R_{k,l},
\end{equation}
where
$$R_{k,l}=\sum_{\underline{k}\in I_{k,l}}\left(1-\frac{1}{k_{m_n+1}!k_{m_n+2}!\dots k_{n-1}!}\right)\prod_{i=m_n+1}^{n-1}\left(1-\frac{i}{n}\right)^{k_i}.$$
Our aim is to determine the first order term of the error when we approximate $S_k$ by $\frac{\lambda_n^k}{k!}$. Since $R_{k,l}\leq S_{k,l}$ for each $l=1,\ldots,k-1$, and for the latter expressions we have the bounds of (\ref{s_maradek}), we see that $\sum_{l=1}^{k-1}R_{k,l}=O\left(\frac{1}{\sqrt{n}}\right)$, and the same, but more detailed argument also gives
\begin{equation}\label{maradek}
\sum_{l=1}^{k-2}R_{k,l}\leq\sum_{l=1}^{k-2}S_{k,l}\leq k!\min\left\{1,(2\lambda_n)^{k}\right\}\frac{1}{n}.
\end{equation}
Thus the leading term of the error $\left|S_k-\frac{\lambda_n^k}{k!}\right|$ is of order $\frac{1}{\sqrt{n}}$, and it comes from the term $R_{k,k-1}$.

Before examining $R_{k,k-1}$ we introduce some notations for further use. As an analogue of the set $I_{k,l}$ we define $I_{k-2,l}$ to be the set of vectors $\underline{k}\in Z_+^{n-m_n-1}$ such that $\sum_{i=m_n+1}^{n-1}k_i=k-2$ and $\underline{k}$ has exactly $l$ nonzero components, $l=1,\ldots,k-2$. Also, as an analogue of the expressions $S_{k,l}$ and $S_k$ we define $S_{k-2,l}$ and $S_{k-2}$ by the formulas in (\ref{s_nm}) with $k$ replaced by $k-2$. Finally we introduce
$$I_{k-2,k-2}^j=\left\{\underline{k}\in I_{k-2,k-2}: k_j=0\right\},\quad j=m_n+1,\ldots,n.$$

We now return to $R_{k,k-1}$. The corresponding index set $I_{k,k-1}$ contains vectors that have exactly one component equal to 2, $k-2$ components equal to 1, and the rest 0. Thus we have
$$R_{k,k-1}
=\frac{1}{2}\sum_{\underline{k}\in I_{k,k-1}}\prod_{i=m_n+1}^{n-1}\left(1-\frac{i}{n}\right)^{k_i}$$
We can write $R_{k,k-1}$ in another form, if we first sum according to the component of the vectors in $I_{k,k-1}$ which equals 2:
\begin{align*}
R_{k,k-1}
&=\frac{1}{2}\sum_{j=m_n+1}^{n}\left(1-\frac{j}{n}\right)^2
\left(\sum_{\underline{k}\in I_{k-2,k-2}^j}\prod_{i=m_n+1}^{n-1}\left(1-\frac{i}{n}\right)^{k_i}\right)\\
&=\frac{1}{2}\sum_{j=m_n+1}^{n}\left(1-\frac{j}{n}\right)^2
\left[\sum_{\underline{k}\in I_{k-2,k-2}}\prod_{i=m_n+1}^{n-1}\left(1-\frac{i}{n}\right)^{k_i}
-\displaystyle{\sum_{\underline{k}\in I_{k-2,k-2}\backslash I_{k-2,k-2}^j}}\prod_{i=m_n+1}^{n-1}\left(1-\frac{i}{n}\right)^{k_i}\right]
\end{align*}
We recognize $S_{k-2,k-2}$ in the first sum in the brackets, thus we can replace it by the formula in (\ref{s_fotag}) with $k-2$ in the place of $k$. As for the second sum in the brackets, we see that $k_j=1$, so there is a $1-\frac{j}{n}$  factor in each of the products, which we can bring before the brackets. These considerations lead to
\begin{align}
R_{k,k-1}
=&\frac{1}{2}\sum_{j=m_n+1}^{n}\left(1-\frac{j}{n}\right)^2\frac{\lambda_n^{k-2}}{(k-2)!}\nonumber\\
&-\frac{1}{2}\sum_{j=m_n+1}^{n}\left(1-\frac{j}{n}\right)^2\sum_{l=1}^{k-3}\sum_{\underline{k}\in I_{k-2,l}}\frac{1}{k_{m_n+1}!k_{m_n+2}!\dots k_{n-1}!}\prod_{i=m_n+1}^{n-1}\left(1-\frac{i}{n}\right)^{k_i}\nonumber\\
&-\frac{1}{2}\sum_{j=m_n+1}^{n}\left(1-\frac{j}{n}\right)^3
\displaystyle{\sum_{\underline{k}\in I_{k-2,k-2}\backslash I_{k-2,k-2}^j}}\prod_{i=m_n+1,i\neq j}^{n-1}\left(1-\frac{i}{n}\right)^{k_i}\nonumber\\
=&:\frac{1}{2}\sum_{j=m_n+1}^{n}\left(1-\frac{j}{n}\right)^2\frac{\lambda_n^{k-2}}{(k-2)!}-R_{k,k-1}^1-R_{k,k-1}^2\label{R_k,k-1}
\end{align}
Now we bound the last two expressions. First,
\begin{equation}\label{R_k,k-1,1}
0\leq R_{k,k-1}^1
\leq\frac{1}{2}\sum_{j=m_n+1}^{n}\left(1-\frac{j}{n}\right)^2\sum_{l=1}^{k-3}S_{k-2,l}
\leq\frac{\lambda_n^{3/2}(k-2)!\min\left\{1,(2\lambda_n)^{k-2}\right\}}{\sqrt{2}}\frac{1}{n}
\end{equation}
by (\ref{lambda_n,j}) and the second inequality in (\ref{s_maradek}) with $k$ replaced by $k-2$. Next,
\begin{equation}\label{R_k,k-1,2}
0\leq R_{k,k-1}^2
\leq\frac{n-m_n-1}{2n}\sum_{j=m_n+1}^{n}\left(1-\frac{j}{n}\right)^2S_{k-2,k-2}
\leq\frac{2^{k-2}\lambda_n^k}{(k-2)!}\frac{1}{n}
\end{equation}
by (\ref{csillag}), (\ref{lambda_n,j}) and the first inequality in (\ref{s_maradek}) with $k$ replaced by $k-2$ and $l=k-2$.

We conclude that if we write (\ref{R_k,k-1}) into (\ref{s_k}), we obtain
\begin{equation*}
S_k=\frac{\lambda_n^k}{k!}
+\frac{1}{2}\sum_{i=m_n+1}^{n}\left(1-\frac{i}{n}\right)^2\frac{\lambda_n^{k-2}}{(k-2)!}
+R_{k,k-1}^1+R_{k,k-1}^2+\sum_{l=1}^{k-2}R_{k,l},
\end{equation*}
where $R_{k,k-1}^1+R_{k,k-1}^2+\sum_{l=1}^{k-2}R_{k,l}=O\left(\frac{1}{n}\right)$ by (\ref{R_k,k-1,1}), (\ref{R_k,k-1,2}), (\ref{maradek}) and the fact that $\lambda_n\to\lambda$ by (\ref{lambda_n}). Thus
\begin{equation}\label{egy}
S_k=\frac{\lambda_n^k}{k!}
+\frac{1}{2}\sum_{i=m_n+1}^{n}\left(1-\frac{i}{n}\right)^2\frac{\lambda_n^{k-2}}{(k-2)!}.
+O\left(\frac{1}{n}\right)
\end{equation}

\bigskip

Now we return to (\ref{p(w=k)}), and approximate the product $\prod_{i=m_n+1}^{n}\frac{i}{n}$ in it by $e^{-\lambda_n}$. Using the definition of $\lambda_n$ in (\ref{lambda_n}) and the expansion formula of the logarithm function the error of the approximation can be written in the form
\begin{align*}
e^{-\lambda_n}-\prod_{i=m_n+1}^{n}\frac{i}{n}
&=\exp\left\{-\sum_{i=m_n+1}^{n}\left(1-\frac{i}{n}\right)\right\}
-\exp\left\{\sum_{i=m_n+1}^{n}\log\left[1-\left(1-\frac{i}{n}\right)\right]\right\}\nonumber \\
&=e^{-\lambda_n}\left(1-\exp\left\{-\sum_{j=2}^{\infty}\frac{1}{j}\lambda_{n,j}\right\}\right)\\
&=e^{-\lambda_n}\left(\frac{1}{2}\sum_{i=m_n+1}^{n}\left(1-\frac{i}{n}\right)^2+\sum_{j=3}^{\infty}\frac{1}{j}\lambda_{n,j}
-\left[\exp\left\{-\sum_{j=2}^{\infty}\frac{1}{j}\lambda_{n,j}\right\}
-1+\sum_{j=2}^{\infty}\frac{1}{j}\lambda_{n,j}\right]\right),
\end{align*}
where the expressions $\lambda_{n,j}$ are defined as in (\ref{lambda_n,j}). Thus we have
\begin{equation}\label{hiba_e}
e^{-\lambda_n}-\prod_{i=m_n+1}^{n}\frac{i}{n}
=e^{-\lambda_n}\frac{1}{2}\sum_{i=m_n+1}^{n}\left(1-\frac{i}{n}\right)^2+R_n,
\end{equation}
where
$$R_n=e^{-\lambda_n}\left(\sum_{j=3}^{\infty}\frac{1}{j}\lambda_{n,j}
-\left[\exp\left\{-\sum_{j=2}^{\infty}\frac{1}{j}\lambda_{n,j}\right\}
-1+\sum_{j=2}^{\infty}\frac{1}{j}\lambda_{n,j}\right]\right),$$
and we are going to show that $R_n=O\left(\frac{1}{n}\right)$.

We are going to bound the sum in the exponent in $R_n$. Since $\lambda_n\to\lambda$ by (\ref{lambda_n}), there exists a threshold number $n_0$ such that for all $n\geq n_0$ we have $\sqrt{\frac{2\lambda_n}{n}}<\frac{1}{2}$. This with inequality (\ref{lambda_n,j}) yields
\begin{equation}\label{lambda_n,j_becslese}
\sum_{j=j_0}^{\infty}\frac{1}{j}\lambda_{n,j}
\leq\lambda_n\left(\frac{2\lambda_n}{n}\right)^{\frac{j_0-1}{2}}\sum_{j=j_0}^{\infty}\left(\sqrt{\frac{2\lambda_n}{n}}\right)^{j-j_0}
\leq\lambda_n\left(\frac{2\lambda_n}{n}\right)^{\frac{j_0-1}{2}}\sum_{j=j_0}^{\infty}\left(\frac{1}{2}\right)^{j-j_0}
=2\lambda_n\left(\frac{2\lambda_n}{n}\right)^{\frac{j_0-1}{2}}
\end{equation}
for all $n\geq n_0$. Let us suppose that $n$ satisfies this condition from now on.

Now we bound $|R_n|$. First we apply the triangle inequality, then the inequality $|e^{-x}-1+x|\leq\frac{x^2}{2}$ valid for all positive real $x$ with $x=\sum_{j=2}^{\infty}\frac{1}{j}\lambda_{n,j}$ , and finally use inequality (\ref{lambda_n,j_becslese}) with $j_0=2$ and 3. Thus we obtain
\begin{align*}
|R_n|
&\leq e^{-\lambda_n}\left(\left|\sum_{j=3}^{\infty}\frac{1}{j}\lambda_{n,j}\right|
+\left|\exp\left\{-\sum_{j=2}^{\infty}\frac{1}{j}\lambda_{n,j}\right\}
-1+\sum_{j=2}^{\infty}\frac{1}{j}\lambda_{n,j}\right|\right)\\
&\leq e^{-\lambda_n}\left(\sum_{j=3}^{\infty}\frac{1}{j}\lambda_{n,j}
+\frac{1}{2}\left(\sum_{j=2}^{\infty}\frac{1}{j}\lambda_{n,j}\right)^2\right)\\
&\leq e^{-\lambda_n}\left(\frac{4\lambda_n^2}{n}
+\frac{1}{2}\left(2\lambda_n\sqrt{\frac{2\lambda_n}{n}}\right)^2\right)=e^{-\lambda_n}4\lambda_n^2(\lambda_n+1)\frac{1}{n}.
\end{align*}
Recalling (\ref{hiba_e}) we see that we proved
\begin{equation}\label{ketto}
e^{-\lambda_n}-\prod_{i=m_n+1}^{n}\frac{i}{n}
=e^{-\lambda_n}\frac{1}{2}\sum_{i=m_n+1}^{n}\left(1-\frac{i}{n}\right)^2+O\left(\frac{1}{n}\right).
\end{equation}

\bigskip

Finally, recalling (\ref{p(w=k)}) we have
$$\p(\widetilde{W}_{n,m_n}=0)=\left(\prod_{i=m_n+1}^{n}\frac{i}{n}\right)=e^{-\lambda_n}-\left(e^{-\lambda_n}-\prod_{i=m_n+1}^{n}\frac{i}{n}\right)$$
for $k=0$,
$$\p(\widetilde{W}_{n,m_n}=1)=\left(\prod_{i=m_n+1}^{n}\frac{i}{n}\right)\lambda_n=e^{-\lambda_n}\lambda_n
-\left(e^{-\lambda_n}-\prod_{i=m_n+1}^{n}\frac{i}{n}\right)\lambda_n$$
for $k=1$, and
$$\p(\widetilde{W}_{n,m_n}=k)=\left(\prod_{i=m_n+1}^{n}\frac{i}{n}\right)S_k
=e^{-\lambda_n}S_k-\left(e^{-\lambda_n}-\prod_{i=m_n+1}^{n}\frac{i}{n}\right)S_k$$
for $k\geq2$. We obtain the first assertion of the Theorem if we write (\ref{egy}) and (\ref{ketto}) into these expressions. The second assertion follows from the first and (\ref{lambda_n}). $\Box$


\begin{thebibliography}{100}


\bibitem{P} P\'osfai, A., Poisson Approximation in a Poisson Limit Theorem Inspired by Coupon Collecting, {\it Journal of Applied Probability} (to appear)

\bibitem{BB} Baum, L. E. and Billingsley, P., Asymptotic distributions for the coupon collector's problem, {\it Ann. Math. Statist.} {\bf 36} (1965), 1835--1839.







\end{thebibliography}
\end{document}